\documentclass[11pt]{article}

\usepackage{changepage} 
\usepackage[colorlinks=true, a4paper=true, pdfstartview=FitV,
linkcolor=blue, citecolor=blue, urlcolor=blue]{hyperref}

\usepackage[normalem]{ulem}
\usepackage{graphicx, amsmath, amssymb, amsthm, float, mathtools}

\usepackage{setspace}

\newtheorem{theorem}{Theorem}
\newtheorem{proposition}{Proposition}
\newtheorem{corollary}{Corollary}
\newtheorem{lemma}{Lemma}

\theoremstyle{definition}
\newtheorem{remark}{Remark}
\newtheorem{example}{Example}

\newcommand{\bbR}[0]{\mathbb{R}}
\newcommand{\bbZ}[0]{\mathbb{Z}}
\newcommand{\transp}[0]{^{\nobreak\hspace{.16em plus .08em}
\intercal}}
\newcommand{\eqt}[1]{\text{\normalfont #1}}
\newcommand{\st}[0]{~ | ~}
\newcommand{\subjectto}[0]{{\normalfont \text{s.t.}}~~}

\newcommand{\ifr}[1]{\mathit{IFR}(#1)}
\newcommand{\inv}[1]{\mathit{Inv}(#1)}
\newcommand{\obj}[1]{z(#1)}
\newcommand{\modw}[0]{(\eqt{mod}~w)}
\newcommand{\gcr}[0]{\mathcal{G}}
\newcommand{\gcrb}[0]{\gcr^{B}}
\newcommand{\gcrbd}[0]{\gcr^{B,d}}
\newcommand{\spp}[0]{\mathcal{S}}
\newcommand{\sppb}[0]{\spp^{B}}
\newcommand{\sppbd}[0]{\spp^{B,d}}
\newcommand{\ip}[0]{\mathcal{IP}}
\newcommand{\ipd}[0]{\ip^d}
\newcommand{\lp}[0]{\mathcal{LP}}
\newcommand{\lpd}[0]{\mathcal{LP}^d}
\newcommand{\ifrgcr}[0]{\ifr{\gcrb, x^\circ}}
\newcommand{\invgcr}[0]{\inv{\gcrb, x^\circ}}
\newcommand{\ifrspp}[0]{\ifr{\sppb, x^\circ_N}}
\newcommand{\invspp}[0]{\inv{\sppb, x^\circ_N}}
\newcommand{\ifrip}[0]{\ifr{\ip, x^\circ}}
\newcommand{\invip}[0]{\inv{\ip, x^\circ}}
\newcommand{\ifrlp}[0]{\ifr{\lp, x^\circ}}
\newcommand{\invlp}[0]{\inv{\lp, x^\circ}}

\usepackage{xcolor, csquotes}

\newcommand{\newv}[1]{#1}

\usepackage[
backend=biber,
style=ieee,
sorting=nyt
]{biblatex}
\DeclareNameAlias{default}{family-given}
\begin{document}

\title{Inverse of the Gomory Corner Relaxation of Integer Programs}
\author{ George Lyu \and Fatemeh Nosrat$^{*}$ \and  Andrew J. Schaefer}

\date{July 2025\\
{\textit{\small{Computational Applied Mathematics and Operations Research, 6100 Main St, Houston, TX 77005, USA}}}\\
\small{ *Corresponding author. E-mail: {\href{mailto:fatemeh.nosrat@rice.edu}{fatemeh.nosrat@rice.edu}}\\
Contributing authors: {\href{mailto:gml8@rice.edu}{gml8@rice.edu}}, {\href{mailto:andrew.schaefer@rice.edu}{andrew.schaefer@rice.edu}}}}

\maketitle
\begin{abstract}
\newv{\noindent We explore the inverse of integer programs (IPs) by studying the inverse of their Gomory corner relaxations (GCRs). We show that solving a set of inverse GCR problems always yields an upper bound on the optimal value of the inverse IP that is at least as tight as the optimal value of the inverse of the linear program (LP) relaxation. We provide conditions under which solving a set of inverse GCR problems exactly solves the inverse IP. We propose an LP formulation for solving the inverse GCR under the $L_1$ and $L_\infty$ norms by reformulating the inverse GCR as the inverse of a shortest path problem.
}

\end{abstract}
\section{Introduction} \label{sec:intro}
Given a (forward) optimization problem and a feasible solution, the \emph{inverse-feasible region} is the set of objective vectors under which the given feasible solution is optimal \newv{for} the forward problem. The \emph{inverse optimization problem} finds an inverse-feasible vector that is closest to a given target \newv{objective} vector \newv{(by some given metric)}. Inverse optimization \newv{applications include healthcare treatment \cite{Ajayi2022, Chan2014, Erkin2010}, medical imaging \cite{Bertero2007}, geophysics \cite{Tarantola2005}, statistics \cite{Tarantola2005}, and traffic modeling \cite{Bertsimas2015}}. 

The inverse of integer programs (IPs) and the inverse of mixed integer programs (MIPs) have been studied widely. Schaefer \cite{Schaefer2009} and Lamperski and Schaefer \cite{Lamperski2015} established polyhedral representations of the inverse-feasible regions of IPs and MIPs using the superadditive duality of the forward problems. This characterization enabled linear program\newv{s} (LP\newv{s}) for inverse IPs and inverse MIPs. However, the number of variables and constraints in these LP formulations grow super-exponentially \newv{on} the size of the forward \newv{IPs/MIPs}. Huang \cite{Huang2005} reformulated the inverse IP as the inverse of a shortest path problem; the number of vertices and arcs in the graph of this shortest path problem grow super-exponentially on the number of constraints in the forward IP.

Cutting plane algorithms have \newv{also} been proposed for solving inverse IPs and MIPs. Wang \cite{Wang2008} provided a cutting plane algorithm for solving inverse MIPs by repeatedly generating optimality cuts from the extreme points of the convex hull of the feasible region of the forward problem. The algorithm was improved upon by Duan and Wang \cite{Duan2011}, who introduced a heuristic algorithm for computing the extreme points and bounds for Wang's algorithm \cite{Wang2008}. Bodur et al. \cite{Bodur2022} introduced another cutting plane algorithm for solving inverse MIPs, which generates optimality cuts from interior points of the convex hull of the feasible region of the forward problem. These cutting plane algorithms are far more tractable than the LP formulations proposed by Schaefer \cite{Schaefer2009} and Lamperski and Schaefer \cite{Lamperski2015}, but the cutting plane algorithms do not characterize the polyhedral structure of the inverse-feasible regions of IPs and MIPs.

IPs and MIPs are often studied by relaxing the integrality constraints, obtaining the LP relaxation. Therefore, a common approach to studying inverse IPs and inverse MIPs is to solve inverse LPs, which typically exhibit more structure. Zhang and Liu \cite{Zhang1996} proposed a solution for general inverse LPs under the $L_1$ norm, from which they obtained strongly polynomial algorithms for solving inverse minimum cost flow problem\newv{s} and inverse assignment problem\newv{s}. Zhang and Liu \cite{Zhang1999} proposed a solution for inverse LPs when both the given feasible solution and an optimal solution under the original objective vector are composed of only zeros and ones, which is common in network flow problems. Ahuja and Orlin \cite{Ahuja2001} showed that if a problem with a linear objective function is polynomially solvable, as is the case for LPs, then the inverse of that problem under the $L_1$ or $L_\infty$ norm is also polynomially solvable. Tavaslıoğlu et al. \cite{Tavaslioglu2018} studied the polyhedral structure of the inverse-feasible region of LPs, while Chan et al. \cite{Chan2019} introduced a goodness-of-fit framework for evaluating inverse LPs where the provided feasible solution for the forward LP cannot be made optimal (outside of the trivial zero-objective case). 

\newv{Given a feasible basis of the LP relaxation, the Gomory corner relaxation (GCR) of the IP with respect to that basis is} obtained by relaxing the nonnegativity constraint of each variable in \newv{that basis} while preserving variable integrality \cite{Gomory1969}. \newv{When the given basis is optimal for the LP relaxation}, Gomory \cite{Gomory1969}, Ho{\c{s}}ten and Thomas \cite{Hosten2003}, and Richard and Dey \cite{Richard2010} enumerated several classes of IP instances where \newv{an} optimal solution for the GCR \newv{is} also \newv{an} optimal solution for the original IP. Fischetti and Monaci \cite{Fischetti2008} demonstrated that for many instances, the gap between the IP and GCR optimal values is much tighter than the gap between the IP and LP relaxation optimal values. Köppe et al. \cite{Koppe2004} characterized the geometry of several reformulations of the GCR. The GCR can be further relaxed to obtain the master group relaxation, \newv{whose general structure} can be applied to broader classes of problems \cite{Richard2010}. The GCR\newv{, when defined on an optimal basis,} is NP-hard \cite{Letchford2003}, and the most efficient known algorithms for solving the GCR exhibit polynomial runtime complexity with respect to the \newv{absolute value of the} determinant of the basis matrix of the LP relaxation, which can be very large \cite{Letchford2003, Richard2010}. Several algorithms for solving the GCR reduce the GCR to an instance of the shortest path problem \cite{Chen1976, Koppe2004, Richard2010}, a technique first developed by Shapiro \cite{Shapiro1968}.

We show that the inverse GCR can be solved as \newv{an inverse shortest path problem}, which manipulates a graph's arc weights such that a given path becomes shortest from among all paths that connect the associated \newv{source} and destination vertices. Inverse shortest path problem\newv{s have} been extensively studied. Forward shortest path problem\newv{s} can be reduced to minimum cost flow problem\newv{s}, so inverse shortest path problem\newv{s} under the $L_1$ norm can be solved using a strongly polynomial algorithm provided by Zhang and Liu \cite{Zhang1996}. Ahuja and Orlin \cite{Ahuja2001} showed that inverse shortest path problem\newv{s} under the $L_1$ norm can be reduced to forward shortest path problem\newv{s}. Burton and Toint \cite{Burton1992} proposed a quadratic programming formulation for solving inverse shortest path problem\newv{s} under the $L_2$ norm. Xu and Zhang \cite{Xu1995} characterized the feasible region of inverse shortest path problem\newv{s} as polyhedral cone\newv{s}.

\newv{An IP may have multiple associated GCR problems, each defined by a different feasible basis of the LP relaxation. We study the inverse IP by solving the inverse GCR for all these feasible bases. Our main contributions are as follows:
\begin{itemize}
    \item In \S \ref{sec:igcr_as_tight}, we prove that solving the inverse GCR for \emph{all} feasible bases provides an upper bound on the optimal value of the inverse IP that is at least as tight as the optimal value of the inverse LP relaxation.
    \item In \S \ref{sec:igcr_as_tight}, we establish mild conditions under which solving the inverse GCR for only \emph{one} feasible basis is sufficient to provide an upper bound on the optimal value of the inverse IP that is at least as tight as the optimal value of the inverse LP relaxation.
    \item In \S \ref{sec:igcr_exact}, we provide conditions under which solving the inverse GCR for \emph{all} feasible bases yields an exact optimal solution for the inverse IP.
    \item In \S \ref{sec:ifr_equiv}, we show that any inverse GCR problem can be reformulated as an inverse shortest path problem.
    \item In \S \ref{sec:igcr_lp}, we introduce an LP formulation for solving any inverse GCR problem. Our formulation has many orders of magnitude fewer variables and constraints than the LP formulation for the inverse GCR obtained by applying the general inverse IP techniques proposed by Schaefer \cite{Schaefer2009}.
\end{itemize}}

\section{\newv{Gomory Corner Relaxation as a Shortest Path Problem}}
\subsection{\newv{Defining the} Gomory Corner Relaxation}\label{sec:define}

Given $A \in \bbZ^{m \times n}$, $b \in \bbZ^m$, and $c \in \bbR^n$, let $\mathcal{IP}$ denote the following IP, which we assume has nonempty feasible region. Let $\mathcal{LP}$ denote the LP relaxation of $\mathcal{IP}$\newv{, and let $\obj{\cdot}$ denote their respective value functions:
\begin{align*}
    \obj{\ip} &:= \min \lbrace c\transp x \st Ax = b, x \geq \mathbf{0}, x \in \bbZ^n \rbrace, \\
    \obj{\lp} &:= \min \lbrace c\transp x \st Ax = b, x \geq \mathbf{0}, x \in \bbR^n \rbrace.
\end{align*}
}Let $B, N \subseteq \{ 1,...,n \}$\ respectively denote the indices of the basic and nonbasic variables of a \newv{given} basic solution for $\mathcal{LP}$. Assume $A$ is full row rank and let $m \leq n$, so $|B| = m, |N| = n - m$ and $B \cap N = \varnothing$. Let $c_B, x_B$ $(c_N, x_N)$ denote the vectors comprised of the $B$-indexed ($N$-indexed) components of $c, x$, respectively. Let $A_B$ ($A_N$) be the matrix comprised of the $B$-indexed ($N$-indexed) columns of $A$. Observe that $A_B$ is nonsingular. Then, the GCR of $\mathcal{IP}$ with respect to $B$, denoted by $\mathcal{G}^B$, is obtained by relaxing the nonnegativity constraints of the decision variables in the selected basis $B$ \cite{Richard2010}:
\newv{\begin{align} 
   \obj{\gcrb} := \min \lbrace c_B\transp x_B + c_N\transp x_N\st A_B x_B + A_N x_N = b, x_N \geq \mathbf{0}, x \in \bbZ^n \rbrace. \label{eq:gcr}
\end{align}}

\subsection{\newv{Shortest Path Reformulation of the Gomory Corner Relaxation}}\label{sec:gcr_reformulation}

\newv{
We reformulate $\gcrb$ as a shortest path problem based on the more comprehensive treatment of Richard and Dey \cite{Richard2010}. We use this shortest path reformulation in \S \ref{sec:igcr_as_tight} to relate the inverse of $\gcrb$ with the inverse of $\lp$ and again in \S \ref{sec:ifr_equiv} to relate the inverse of $\gcrb$ with the inverse of the shortest path problem.

\begin{lemma} \label{lem:snf} {\normalfont \textbf{\cite{Richard2010}} 
    (Smith Normal Form)} For any nonzero matrix $A_B \in \bbZ^{m \times m}$, there exists a positive vector $w \in \bbZ^m_{++}$ and unimodular matrices $S, T \in \bbZ^{m \times m}$ such that $S A_B T = \text{\normalfont diag}(w)$, where $\text{\normalfont diag}(w)$ is the $m \times m$ matrix whose diagonal is given by $w$ and whose off-diagonal entries are all zero.
\end{lemma}
 
Lemma \ref{lem:snf} provides the Smith Normal Form of $A_B$ \cite{Smith1861}, for which there exists several efficient algorithms \cite{Birmpilis2023, Dumas2001}. $S, T, w$ depend on the selected basis $B$, but our notation suppresses this dependence for clarity. Observe that for any feasible solution $x$ for $\gcrb$, we have $x_B = A_B^{-1} (b - A_N x_N)$. Then,
\begin{align}
       \obj{\gcrb} &= \min \lbrace c_B\transp A_B^{-1} (b - A_N x_N)  + c_N\transp x_N\st S A_B x_B + S A_N x_N = S b, x_N \geq \mathbf{0}, x \in \bbZ^n \rbrace \notag \\
       &= \min \lbrace c_B\transp A_B^{-1} b + (c_N\transp  -c_B\transp A_B^{-1} A_N) x_N \st \text{\normalfont diag}(w) T^{-1} x_B + S A_N x_N = S b, x_N \geq \mathbf{0}, x \in \bbZ^n \rbrace. \notag
\end{align}

$T^{-1}$ is unimodular because $T$ is unimodular, and $x_B$ has all integer entries, so $T^{-1} x_B$ has all integer entries. Thus, the $j$th entry of the vector $\text{\normalfont diag}(w) T^{-1} x_B$ must be an integer multiple of $w_j$ for all $j \in [n-m]$. Additionally, let $\bar c_N := (c_N\transp  -c_B\transp A_B^{-1} A_N)^\intercal$ denote the reduced costs of the $N$-indexed variables in the basic solution at basis $B$ of $\lp$. Then, $\gcrb$ can be rewritten as
\begin{align}
    \obj{\gcrb} = \min \lbrace \bar c_N \transp x_N \st S A_N x_N \equiv S b \modw, x_N \in \bbZ^{n-m}_+ \rbrace + c_B\transp A_B^{-1} b, \label{eq:gcr_reform}
\end{align}
where for arbitrary given $u, v \in \bbZ^{n-m}$, the mod notation $u \equiv v \modw$ denotes that $u_j - v_j$ is an integer multiple of $w_j$ for all $j \in [n-m]$ \cite{Richard2010}.

Under $\modw$ equivalence, there are at most $\prod_{j \in [n-m]} w_j$ distinct values that the left-hand side $S A_N x_N$ can take in \eqref{eq:gcr_reform}. That is, for any value of $S A_N x_N$, we can produce a vector $v \equiv S A_N x_N \modw$ with $v_j \in \{0, ..., w_j-1\}$ for each $j \in [n-m]$ by setting $v_j$ to be $(S A_N x_N)_j$ modulo $w_j$. Let $V$ be the set of all distinct elements that exist under $\modw$ equivalence; i.e., using the Cartesian product, $V := \prod_{j \in [n-m]} \{0, ..., w_j - 1\}$. Consider the directed graph with vertex set $V$ and arc set $E := \bigcup_{j \in [n-m]} E_j$, where $E_j := \lbrace (u, u + (S A_N)_j) \st u \in V \rbrace$ and $(S A_N)_j$ is the $j$th column of the matrix $S A_N$.

Define $\sppb$ as the shortest path problem from the source vertex $\mathbf{0}$ to the destination vertex $Sb$, where each arc from $E_j$ has weight $(\bar{c}_N)_j$. We represent solutions to $\sppb$ as vectors from $\bbZ_+^{n-m}$, where the vector $p \in \bbZ_+^{n-m}$ represents a path starting at source vertex $\mathbf{0}$ that traverses exactly $p_j$ arcs from $E_j$ for each $j \in [n-m]$. Such a path always exists because each vertex is the tail vertex of exactly one arc from $E_j$ for each $j$. The vector $p$ may represent multiple different paths because $p$ does not encode the order in which the arcs are traversed, but these different paths all have the same destination vertex $S A_N p$ and the same total weight $\bar c_N \transp p$, so these different paths all have the same feasibility and objective value for $\sppb$.

Richard and Dey \cite{Richard2010} showed that a solution $x_N$ for \eqref{eq:gcr_reform} is analogous to a path $x_N$ for $\sppb$, as captured by Lemmas \ref{lem:gcr_equiv_sp} and \ref{lem:gcr_objval_sp}, when $B$ is an optimal basis of $\lp$. A similar proof holds when $B$ is nonoptimal: First, $S A_N x_N \equiv S b \modw$ if and only if the path $x_N$ terminates at vertex $S b$ for $\sppb$, so $x_N$ is feasible for \eqref{eq:gcr_reform} if and only if $x_N$ is a feasible $\mathbf{0}$-to-$S b$ path for $\sppb$. Second, $x_N$ has the objective value $\bar c_N \transp x_N + c_B\transp A_B^{-1} b$ in \eqref{eq:gcr_reform}, which differs from the weight $\bar c_N \transp x_N$ of path $x_N$ for $\sppb$ by exactly the fixed difference $c_B\transp A_B^{-1} b$.

\begin{lemma} \label{lem:gcr_equiv_sp} {\normalfont\textbf{\cite{Richard2010}}}
    Consider any $x \in \bbZ^n$. Then, $x$ is an optimal solution for $\gcrb$ if and only if $x_N$ is an optimal solution for $\sppb$ and $x_B = A_B^{-1} (b - A_N x_N)$.
\end{lemma}

\begin{lemma} \label{lem:gcr_objval_sp}{\normalfont\textbf{\cite{Richard2010}}}
    We have $\obj{\gcrb} = \obj{\sppb} + c_B\transp A_B^{-1} b$. If $B$ is an optimal basis of $\lp$, we have $\obj{\gcrb} = \obj{\sppb} + \obj{\lp}$ because $\obj{\lp} = c_B\transp A_B^{-1} b$.    
\end{lemma}

Since $S$ and $T$ are unimodular, $|V| = \prod_{j \in [n-m]} w_j = |\det A_B|$ and $|E| = (n-m) |V| = (n-m) |\det A_B|$. We will use this shortest path reformulation $\sppb$ of $\gcrb$ to reformulate the inverse of $\gcrb$ as the inverse of $\sppb$ in \S \ref{sec:ifr_equiv}.
    
\begin{example} \label{ex:forward}

Consider the following IP instance:

\begin{equation*}
\begin{array}{llllll}
    \min \quad & & & -2x_3& - 3x_4& \\
    \subjectto & x_1 & \quad & + 2x_3& + 4x_4 & = 9, \\
    & & x_2 & + 4x_3 & + 4x_4& = 15, \\
    & & & x \in \bbZ_+^n.
\end{array}    
\end{equation*}
We select basis $B = \{3, 4\}$, $N = \{1, 2\}$. Then, $\gcrb$ relaxes the nonnegativity constraint of $x_3, x_4$ (Figure \ref{fig:example}). It is easy to verify that $\obj{\ip} = \obj{\gcrb} = -7$ with optimal solution $(1, 3, 2, 1)\transp$, and $\obj{\lp} = -33/4$ with optimal solution $(0, 0, 3, 3/4)\transp$. $\gcrb$ solves $\ip$ exactly, whereas $\lp$ does not.

\begin{figure}[H]
    \centering
    \includegraphics[height=6.5cm]{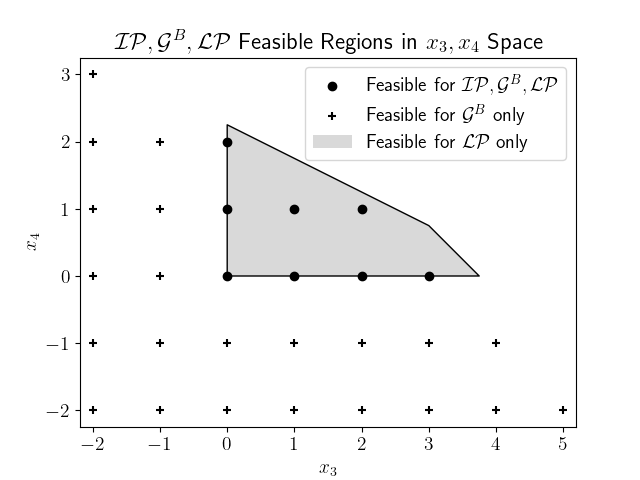}
    \includegraphics[height=6.5cm]{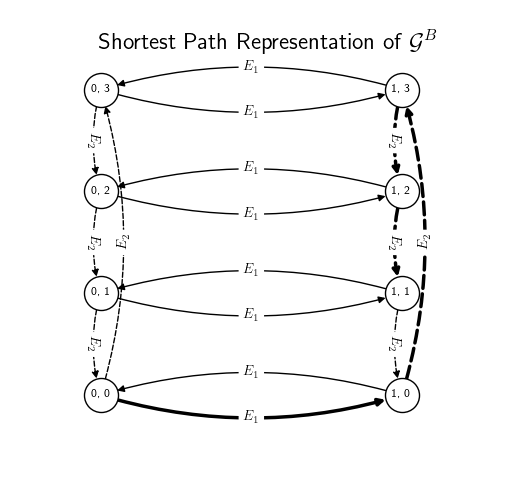}
    \caption{\newv{\emph{(Left)} Feasible regions of $\ip, \gcrb, \lp$ from Example \ref{ex:forward} are projected onto the $x_3, x_4$ space. Feasible solutions for $\ip$ are shown using round markers. Feasible solutions for $\gcrb$ that are infeasible for $\ip$ are shown using plus markers. The feasible region of $\lp$ is shaded. \emph{(Right)} The directed graph shows the reformulation of $\gcrb$ as $\sppb$ in Example \ref{ex:forward}. Arcs from $E_1$ ($E_2$) are shown using solid (dashed) arrows. The shortest paths from $\mathbf{0}$ to $S b \equiv (1, 1) \transp \modw$ traverse one arc from $E_1$ and three arcs from $E_2$. One such shortest path is shown in bolded arcs.}}
    \label{fig:example}
\end{figure}

The Smith Normal Form of $A_B = \begin{pmatrix} 2 & 4 \\ 4 & 4 \end{pmatrix}$ is given by
$    S = \begin{pmatrix}
        1 & 0 \\
        0 & -1
    \end{pmatrix}, ~~~
    T = \begin{pmatrix}
        -1 & -2 \\
        1 & 1
    \end{pmatrix}, ~~~
    w = \begin{pmatrix}
        2 \\
        4
    \end{pmatrix}.
$ Then, $S A_N = \begin{pmatrix} 1 & 0 \\ 0 & -1 \end{pmatrix}$ and $S b = \begin{pmatrix} 9 \\ -15 \end{pmatrix} \equiv \begin{pmatrix} 1 \\ 1 \end{pmatrix} \modw$. We have $\bar c_N = \begin{pmatrix} 1/2 \\ 1/4 \end{pmatrix}$.

Thus, $\sppb$ is problem of finding the shortest path from vertex $\mathbf{0}$ to vertex $(1, 1)\transp$ on the directed graph shown in Figure \ref{fig:example}. The arc from $E_1$ that starts at arbitrary vertex $u \in V$ will terminate at vertex $u + (SA_N)_1 = (u_1 + 1, u_2)\transp$ and incur a weight of $(\bar c_N)_1 = 1/2$. Likewise, the arc from $E_2$ that starts at vertex $u$ will terminate at vertex $u + (SA_N)_2 = (u_1, u_2 - 1)\transp$ and incur a weight of $(\bar c_N)_2 = 1/4$. The shortest paths traverse one arc from $E_1$ and three arcs from $E_2$, which are exactly the $N$-indexed components of the optimal solution for $\ip$ and $\gcrb$.
\end{example}}

\newv{
\section{Tighter Upper Bounds on Inverse IP Optimal Value}\label{sec:inverse}
}

\newv{For any arbitrary $d \in \bbR^n$, let $\ipd, \lpd, \gcrbd$ denote the problems $\ip, \lp, \gcrb$, respectively, after replacing the objective vector with $d$. Let $\sppbd$ denote $\sppb$ after replacing the weight of each arc from $E_j$ with $(\bar d_N)_j$ for each $j \in [n-m]$.
}

Let $\mathcal{P}$ be an optimization problem from among $\mathcal{IP}, \mathcal{LP}, \mathcal{G}^B$. Let $x^\circ$ be a feasible solution for $\mathcal{P}$. The \newv{\emph{inverse-feasible region}} of $\mathcal{P}$ with respect to $x^\circ$, denoted by $\ifr{\mathcal{P}, x^\circ}$, is the set of vectors $d \in \bbR^n$ for which $x^\circ$ is an optimal solution for $\mathcal{P}^{d}$:
\begin{align*}
    \ifr{\mathcal{P}, x^\circ} ~\newv{:=} ~\lbrace d \in\mathbb{R}^{n} \st \newv{d \transp x^\circ = \obj{\mathcal{P}^d}} \rbrace.
\end{align*} The \newv{\emph{inverse problem}} of $\mathcal{P}$ with respect to $x^\circ$, denoted by $\inv{\mathcal{P}, x^\circ}$, is the problem of finding a vector $d \in \ifr{\mathcal{P}, x^\circ}$ that minimizes the (possibly weighted) $L_p$ \newv{distance between $d$ and the \emph{target objective vector} $c$}:
\newv{\begin{align*}
\obj{\inv{\mathcal{P}, x^\circ}} := \min\limits \lbrace \lVert d - c \rVert_p \st d \in \ifr{\mathcal{P}, x^\circ} \rbrace.
\end{align*}}\newv{We similarly define the inverse-feasible region of $\sppb$ and the inverse of $\sppb$, which require slight modification because the feasible solutions of $\sppb$ have $n-m$ components (whereas the feasible solutions of $\ip, \lp, \gcrb$ have $n$ components). Let $x^\circ_N$ be a feasible $\mathbf{0}$-to-$Sb$ path for $\sppb$. Then, we define the inverse-feasible region and inverse problem of $\sppb$ with respect to $x^\circ_N$ as
\begin{align*}
    \ifrspp &:= \lbrace d \in \bbR^n \st \bar d_N\transp x_N^\circ = \obj{\sppbd} \rbrace, \\
    \obj{\invspp} &:= \min\limits \lbrace \lVert d - c \rVert_p \st d \in \ifrspp \rbrace.
\end{align*}

\subsection{Inverse Gomory Corner Relaxation Bounds are at Least as Tight as Those of the Inverse LP Relaxation} \label{sec:igcr_as_tight}

We compare the inverses of $\gcrb$ and $\lp$. Let $\mathcal{B}$ be the (finite) set of all feasible bases of $\lp$.

\begin{proposition} \label{prop:lp_opt_basis}
    Consider any $B \in \mathcal{B}$, any feasible solution $x^\circ$ for $\gcrb$, and any $d \in \ifrgcr$. Then, $B$ is an optimal basis of $\lpd$.
\end{proposition}
\begin{proof} 
    By contradiction, suppose there exists $d \in \ifrgcr$ such that $B$ is a nonoptimal feasible basis of $\lpd$. Nonoptimality implies there exists a negative reduced cost; that is, there exists $j \in [n-m]$ such that $(\bar d_N)_j < 0$ \cite{Bertsimas1997}. Then, $\gcrbd$ is unbounded \cite{Richard2010}, as the vector $y \in \bbZ_+^{n-m}$ defined by $y_j := \prod_{k \in [n-m]} w_k$ and $y_l := 0$ for $l \in [n-m] \setminus \{j\}$ is an unbounded ray for \eqref{eq:gcr_reform} because $S A_N y \equiv \mathbf{0} \modw$ and $\bar d_N \transp y < 0$. Thus, $\gcrbd$ has no optimal solution and $d \not\in \ifrgcr$.
\end{proof}

\begin{remark}
    Our treatment of $\gcrb$ in \S \ref{sec:define} and \S \ref{sec:gcr_reformulation} allows $B$ to be a nonoptimal basis of $\lp$. Many authors focus on the case where $B$ is an optimal basis \cite{Fischetti2008, Richard2010}. Proposition \ref{prop:lp_opt_basis} reconciles these treatments, since $B$ must be an optimal basis of $\lpd$ for all $d \in \ifrgcr$.
\end{remark}

If $x^\circ$ is a feasible solution for $\ip$, then $\invip$ finds an objective vector $d$ such that $x^\circ$ is an optimal solution for $\ipd$. Theorem \ref{thm:igcr_as_tight} uses Proposition \ref{prop:lp_opt_basis} to show that solving all $\invgcr, B \in \mathcal{B}$ provides an upper bound on the optimal value of $\invip$ that is at least as tight as the optimal value of $\invlp$.

\begin{theorem} \label{thm:igcr_as_tight}
    $\ifrlp = \big\lbrace d \in \bigcup_{B \in \mathcal{B}} \ifrgcr \st \bar d_N \transp x^\circ_N = 0 \big \rbrace$, where $x^\circ$ is any given feasible solution for $\ip$.
\end{theorem}
\begin{proof}
    First, we show that $\ifrlp \subseteq \bigcup_{B \in \mathcal{B}} \ifrgcr$. Consider arbitrary $d \in \ifrlp$. Then, $x^\circ$ is an optimal solution for $\lpd$ and there exists an optimal basis $\hat B \in \mathcal{B}$ for $\lpd$ such that the reduced costs of $\lpd$ are nonnegative under this basis; that is, $\bar d_{\hat N} \geq \mathbf{0}$ for $\hat N := [n] \setminus \hat B$ \cite{Richard2010}. Thus, the weights of the arcs for $\spp^{\hat B, d}$ are nonnegative and $\obj{\spp^{\hat B, d}} \geq 0$. By Lemma \ref{lem:gcr_objval_sp}, $\obj{\gcr^{\hat B, d}} = \obj{\spp^{\hat B, d}} + \obj{\lpd} \geq \obj{\lpd} = d\transp x^\circ$. Since $x^\circ$ is a feasible solution for $\gcr^{\hat B, d}$, we have $\obj{\gcr^{\hat B, d}} \leq d \transp x^\circ$. Therefore, $\obj{\gcr^{\hat B, d}} = d \transp x^\circ$ and $d \in \ifr{\gcr^{\hat B, d}, x^\circ}$.

    Next, considering arbitrary $d \in  \bigcup_{B \in \mathcal{B}} \ifrgcr$, we show that $d \in \ifrlp$ if and only if $\bar d_N \transp x^\circ = 0$. There exists a basis $\hat B \in \mathcal{B}$ such that $d \in \ifr{\gcr^{\hat B, d}, x^\circ}$ and $\obj{\gcr^{\hat B, d}} = d \transp x^\circ$. Proposition \ref{prop:lp_opt_basis} shows that $\hat B$ is an optimal basis of $\lpd$. Thus, $d \in \ifrlp$ if and only if $\obj{\lpd} = d \transp x^\circ$ if and only if $\obj{\gcr^{\hat B, d}} = \obj{\lpd}$ if and only if, by Lemma \ref{lem:gcr_objval_sp}, $\obj{\spp^{\hat B, d}} = 0$ if and only if $\bar d_{\hat N} \transp x^\circ_{\hat N} = 0$.
\end{proof}

\begin{corollary}
    Since $\gcrb$ is a relaxation of $\ip$ for any $B \in \mathcal{B}$, we know $x^\circ$ is an optimal solution for $\ip$ if $x^\circ$ is an optimal solution for $\gcrb$. Thus, $\ifrlp \subseteq \bigcup_{B \in \mathcal{B}} \ifrgcr \subseteq \ifrip$. Additionally, $\obj{\invlp} \geq \min_{B \in \mathcal{B}} \obj{\invgcr} \geq \obj{\invip}$.
\end{corollary}

\begin{corollary}
    $\ifrlp \subsetneq \bigcup_{B \in \mathcal{B}} \ifrgcr$ if and only if there exist $\hat B \in \mathcal{B}$, $\hat N := [n] \setminus \hat B$ and $d \in \ifr{\gcr^{\hat B}, x^\circ}$ such that $\bar d_{\hat N} \transp x^\circ_{\hat N} > 0$.
\end{corollary}

Whereas Theorem \ref{thm:igcr_as_tight} utilizes all feasible bases of $\lp$, Theorem \ref{thm:igcr_basis} identifies a specific basis $B \in \mathcal{B}$ where $\invgcr$ may provide an upper bound on the optimal value of $\invip$ that is at least as tight as the optimal value of $\invlp$.

\begin{theorem} \label{thm:igcr_basis}
    If $x^\circ$ is a feasible solution for $\ip$, and $x^\circ_B > \mathbf{0}$ for a given $B \in \mathcal{B}$, then $\ifrlp \subseteq \ifrgcr$.
\end{theorem}
\begin{proof}  
By contrapositive, consider arbitrary $d \in \bbR^n$ where $d \notin \ifrgcr$. Then, $x^\circ$ is not an optimal solution for $\gcrbd$, so there exists a feasible solution $s$ for $\gcrbd$ such that $d \transp s < d \transp x^\circ$. We will show that $x^\circ$ is not an optimal solution for $\lpd$.

\noindent \emph{Case 1.} Suppose $s \geq \mathbf{0}$. Because $s$ is a feasible solution for $\gcrb$, we have $As = b$, so $s$ is a feasible solution for $\lpd$. Then, $x^\circ$ is not optimal for $\lpd$ because $d \transp s < d \transp x^\circ$.
    
\noindent \emph{Case 2.} Let $C := \{k \in B \st s_k < 0\} \neq \varnothing$ denote the indices of the negative components of $s$. Let $\lambda := \max \{ -s_k/(x^\circ_k - s_k) ~ | ~ k \in C \}$. We have $\lambda \in (0, 1)$ because $s_k < 0 < x^\circ_k$ for each $k \in C$. Thus, $w := s + \lambda(x^\circ - s) \geq \mathbf{0}$ is a linear combination of $x^\circ, s$. Additionally $Aw = b$ because $Ax = As = b$, so $w$ is a feasible solution for $\lp$. Furthermore, $\lambda \neq 1$ and $d \transp s < d \transp x^\circ$; therefore, $d \transp w < d \transp x^\circ$, and $x^\circ$ is not an optimal solution for $\lpd$.
\end{proof}

Theorem \ref{thm:igcr_basis} shows that if $x^\circ_B$ has all positive components for some basis $B \in \mathcal{B}$, then $\ifrlp \subseteq \ifrgcr$ and $\obj{\invlp} \geq \obj{\invgcr} \geq \obj{\invip}$. Note that $|B| = m$, and such a basis must exist if $x^\circ$ is a nondegenerate basic feasible solution for $\lp$.
\\
\\
\noindent\textbf{Example \ref{ex:forward} Revisited.} Let $x^\circ = (1, 3, 2, 1)\transp$, which is an optimal solution for both $\ip$ and $\gcrb$ for $B = \{3, 4\}$, so $c \in \ifrip, \ifrgcr$. Observe that $x^\circ$ is an interior point of the feasible region of $\lp$, so $c \not \in \ifrlp$. Also observe that $x^\circ_B > \mathbf{0}$; thus, by Theorem \ref{thm:igcr_basis}, $\ifrlp \subsetneq \ifrgcr$.

\subsection{Solving the Inverse IP Using the Inverse Gomory Corner Relaxation} \label{sec:igcr_exact}

Considering a given feasible solution $x^\circ$ for $\ip$, we provide conditions where solving all $\invgcr$, $B \in \mathcal{B}$ exactly solves $\invip$. Gomory \cite{Gomory1969} identified conditions for when the forward problem $\gcrb$ exactly solves $\ip$ for some $B \in \mathcal{B}$, which we paraphrase in Lemma \ref{lem:asymptotic}.

\begin{lemma} \label{lem:asymptotic} {\normalfont\textbf{\cite{Gomory1969}}}
    Suppose $B$ is an optimal basis of $\lp$ and that $\ip$ has a feasible solution. Suppose the $L_2$ distance between the right-hand side $b$ and the boundary of the cone $A_B^{-1} x \geq \mathbf{0}$ is at least $|\det A_B| \max_{j \in [n-m]} \lVert (A_N)_j \rVert_2$. Then, $\obj{\ip} = \obj{\gcrb}$.
\end{lemma}

Thus, we have conditions for when solving all $\invgcr$, $B \in \mathcal{B}$ will exactly solve $\invip$.

\begin{theorem} \label{thm:igcr_asymptotic}
    Suppose $x^\circ$ is a feasible solution for $\ip$, and $\ip$ satisfies the conditions of Lemma \ref{lem:asymptotic}. Then, $\bigcup_{B \in \mathcal{B}} \ifrgcr = \ifrip$.
\end{theorem}
\begin{proof}
    Since $\gcrb$ is a relaxation of $\ip$ for all $B \in \mathcal{B}$, we have $\bigcup_{B \in \mathcal{B}} \ifrgcr \subseteq \ifrip$. To prove equality, consider any $d \in \ifrip$. Since $\ipd$ has an optimal solution $x^\circ$ and $A$ has integer entries, then $\lpd$ must have an optimal solution \cite{Byrd1987}. Let $\hat B$ be an optimal basis of $\lpd$. By assumption, $\ip$ satisfies the conditions of Lemma \ref{lem:asymptotic}; thus, $x^\circ$ must be an optimal solution for $\gcr^{\hat B, d}$. Therefore, $d \in \ifr{\gcr^{\hat B, d}, x^\circ}$ and $\ifr{\gcr^{\hat B, d}, x^\circ} \supseteq \ifrip$.
\end{proof}

\begin{remark}
    When $\bigcup_{B \in \mathcal{B}} \ifrgcr = \ifrip$, then we have $\min_{B \in \mathcal{B}} \obj{\invgcr} = \obj{\invip}$.
\end{remark}
}

\section{\newv{LP Formulation for the Inverse Gomory Corner Relaxation}}\label{sec:inverse_gcr_lp}

\subsection{\newv{Inverse Shortest Path Formulation for the Inverse Gomory Corner Relaxation}} \label{sec:ifr_equiv}

\newv{Consider some basis $B \in \mathcal{B}$ and some feasible solution $x^\circ$ for $\gcrb$. We reformulate $\invgcr$ as an inverse shortest path problem, mirroring how the forward $\gcrb$ can be reformulated as a forward shortest path problem. First, we present a polyhedral representation of $\ifrspp$, where $x_N^\circ$ encodes a $\mathbf{0}$-to-$S b$ path that traverses exactly $(x_N^\circ)_j$ arcs from $E_j$ for each $j \in [n-m]$.} Lemma \ref{lem:ahuja_network} follows from applying known conditions for \newv{$x_N^\circ$ to be a shortest path} for problem $\mathcal{S}^{B,d}$ (e.g., see Chapter 5.2 in Ahuja et al. \cite{Ahuja1993}).

\begin{lemma} \label{lem:ahuja_network}
    For a given $d \in \bbR^n$, \newv{let $x_N^\circ$ be a $\mathbf{0}$-to-$S b$ path for $\sppbd$.} Then, $x_N^\circ$ is a shortest $\mathbf{0}$-to-$S b$ path if and only if for each \newv{graph} vertex $u \in V$, there exists an associated $y_u \in \bbR$ such that
    \newv{
        \begin{subequations}
    \begin{alignat}{2}
        y_{\mathbf{0}} &= 0, \label{eq:isp_source} \\
        y_{S b} &= \bar d_N \transp x_N^\circ, \label{eq:isp_destination} \\
        y_v - y_u &\leq (\bar d_N)_j, & \quad \forall (u, v) \in E_j, ~ \forall j \in [n - m]. \label{eq:isp_dual}
    \end{alignat}
    \end{subequations}
    }
\end{lemma}
$\ifrspp$ is the set of all $d \in \bbR^n$ such that $x^\circ_N$ is a shortest $\mathbf{0}$-to-$S b$ path for problem $\mathcal{S}^{B,d}$, so we formulate $\ifrspp$ by defining the set of all $d \in \bbR^n$ that satisfy the conditions in Lemma \ref{lem:ahuja_network} given by \eqref{eq:isp_source}, \eqref{eq:isp_destination}, \eqref{eq:isp_dual}\newv{, which are linear constraints with respect to $d$ and $y$.}

\begin{proposition} \label{prop:lp_sd}
\newv{
Let $x_N^\circ$ be a $\mathbf{0}$-to-$S b$ path for $\sppb$. Then,
\begin{align*}
    \ifrspp &= \text{\normalfont proj}_{\bbR^n} \lbrace d \in \bbR^n, y \in \bbR^{|V|} \st \eqref{eq:isp_source}, \eqref{eq:isp_destination}, \eqref{eq:isp_dual} \rbrace \\
    &= \lbrace d \in \bbR^n \st \exists y \in \bbR^{|V|} \text{\normalfont ~such that~} \eqref{eq:isp_source}, \eqref{eq:isp_destination}, \eqref{eq:isp_dual} \rbrace.
\end{align*}
}which is a polyhedral cone that contains $\mathbf{0}$.
\end{proposition}

Note that $\invspp$ is a special case of the general inverse shortest path problem (e.g., studied by Zhang and Liu \cite{Zhang1996} and Ahuja and Orlin \cite{Ahuja2001}) because\newv{, for each $j \in [n-m]$, all the} arcs in $E_j$ \newv{must have} the same weight. \newv{We represent $\ifrgcr$ using $\ifrspp$ in Theorem \ref{thm:ifr_equiv}.}

\begin{theorem} \label{thm:ifr_equiv}
For a given feasible solution $x^\circ$ for $\mathcal{G}^B$, $\ifrspp = \ifrgcr$.
\end{theorem}
\begin{proof}
Since $x^\circ$ is a feasible solution for $\mathcal{G}^B$, we have  $x_B^\circ = A_B^{-1} b - A_B^{-1} A_N x_N^\circ$. Then, by Lemma \ref{lem:gcr_equiv_sp}, \newv{$d \in \ifrgcr$ if and only if $x^\circ$ is an optimal solution for $\gcrbd$ if and only if $x^\circ_N$ is an optimal solution for $\sppbd$ if and only if $d \in \ifrspp$.}
\end{proof}

\subsection{\newv{LP Formulation for the Inverse Gomory Corner Relaxation}}\label{sec:igcr_lp}

We obtain the following LP formulation for $\invgcr$ under the $L_1$ norm. The constraints are derived from Proposition \ref{prop:lp_sd}, and we linearize $\lVert d - c \rVert_1$ by substituting $d := c - e + f$ for $e, f \in \bbR_+^n$.

\begin{proposition} \label{prop:lp_formulation}
    For a given feasible solution $x^\circ$ for $\mathcal{G}^B$, an optimal solution for $\invgcr$ under the $L_1$ norm weighted by a given $w \in \bbR_+^n$ is equal to $c - e^* + f^*$, where $e^*, f^*, y^*$ is an optimal solution for the following LP:
    \newv{
    \begin{subequations} \label{eq:igcr_lp}
    \begin{alignat}{3}
    &\mathrm{min} \quad & w \transp (&e + f) \\
    &\subjectto \quad & y_{\mathbf{0}} &= 0, \\ 
    &  & y_{S b} &= (\bar c_N - \bar e_N + \bar f_N) \transp x_N^\circ, \\
    & &y_v - y_u &\leq (\bar c_N)_j - \bar (\bar e_N)_j + (\bar f_N)_j, & \quad \forall(u, v) \in E_j, ~ \forall j \in [n-m], \\
    & &e, f &\in \bbR_+ ^n, y \in \bbR^{|V|}.
    \end{alignat}
    \end{subequations}}
\end{proposition}

The LP \eqref{eq:igcr_lp} can be modified to solve $\invgcr$ under the $L_\infty$ norm \newv{in a manner similar to that of} \cite{Ahuja2001}. \newv{The number of variables and constraints in \eqref{eq:igcr_lp} depends on $|V|$ and $|E|$, which} can be quite large depending on $\det A_B$.

\newv{When $\ip$ satisfies the conditions for Lemma \ref{lem:asymptotic}, then Proposition \ref{prop:lp_formulation} and Theorem \ref{thm:igcr_asymptotic} indicate a finitely terminating algorithm for solving $\invip$. In this case, $\min_{B \in \mathcal{B}} \obj{\invgcr} = \obj{\invip}$, which we can compute by solving all $\invgcr$, $B \in \mathcal{B}$ using the LP in \eqref{eq:igcr_lp}.}

\subsection{\newv{Inverse Gomory Corner Relaxation Formulation Size Benchmark}}
\label{sec:igcr_compare}

\newv{Let $x^\circ$ be a feasible solution for $\gcrb$. Schaefer \cite{Schaefer2009} introduced an LP formulation for exactly solving general inverse IPs. Thus, we can obtain an LP formulation for $\invgcr$ by treating $\gcrb$ as an IP instance and applying Schaefer's \cite{Schaefer2009} LP formulation. We benchmark the number of variables and constraints (``size") in our formulation for $\invgcr$ from \eqref{eq:igcr_lp} by comparing it with Schaefer's \cite{Schaefer2009} formulation for $\invgcr$. Our goal is not necessarily to show that $\invgcr$ yields smaller LP formulations than $\invip$, but rather to show that our LP formulation for $\invgcr$ is significantly smaller than an existing benchmark formulation for $\invgcr$ obtained by applying an LP formulation for general inverse IPs to $\gcrb$.}

\newv{\begin{table}[H]
\caption{\newv{We compare our LP formulation for $\invgcr$ in \eqref{eq:igcr_lp} with Schaefer's \cite{Schaefer2009} for five IP instances from \cite{Miplib}. ``$\ip$" lists the size of the original IP instance. ``$\invgcr$ by \eqref{eq:igcr_lp}" lists the size of our LP formulation for $\invgcr$ obtained by converting the IP into the form of \eqref{eq:gcr} and applying \eqref{eq:igcr_lp}. ``$\invgcr$ by \cite{Schaefer2009}" lists the size of Schaefer's \cite{Schaefer2009} formulation for $\invgcr$ obtained by converting $\gcrb$ from the form of \eqref{eq:gcr} into inequality form and applying his formulation for general inverse IPs.} \\ }
\centering
{\begin{tabular}{|lrr|rr|rr|}
\hline
\multicolumn{3}{|c}{$\ip$} & \multicolumn{2}{|c|}{$\invgcr$ by \eqref{eq:igcr_lp}}      & \multicolumn{2}{c|}{$\invgcr$ by \cite{Schaefer2009}}  \\
Instance       & var  & con  & $\log_{10}$ var & $\log_{10}$ con & $\log_{10}$ var & $\log_{10}$ con \\
\hline
gen-ip016  & 28     & 24     & 2.9           & 4.3           & 105.8         & 197.6         \\
gen-ip054  & 30     & 27      & 11.6          & 13.0          & 77.6          & 141.0       \\
gen-ip021  & 35     & 28      & 10.1          & 11.7          & 104.6         & 193.0     \\
gen-ip002  & 41     & 24      & 20.1         & 21.7          & 103.1         & 192.2    \\
ns1952667  & 13264  & 41      & 32.8          & 36.9          & 244.5         & 464.7       \\
\hline
\end{tabular}}
\label{tbl:compare_lp}
\end{table}}

\newv{We make} this comparison for each of five pure IP instances obtained from MIPLIB 2017 \cite{Miplib}. For each instance, $B$ is an optimal basis of the LP relaxation, computed using Gurobi 10.0.2 \cite{Gurobi2023}. \newv{$x^\circ$ is set to an arbitrary feasible solution for $\gcrb$, as the size of both our formulation in \eqref{eq:igcr_lp} and Schaefer's \cite{Schaefer2009} formulation are independent of $x^\circ$. The goal of $\invgcr$ is to find an objective vector $d$ under which $x^\circ$ is optimal for the forward problem $\gcrbd$.}

Our LP formulation \newv{in \eqref{eq:igcr_lp}} has $2n + |\det A_B|$ variables and $2 + (n-m) |\det A_B|$ constraints. Schaefer's \cite{Schaefer2009} LP formulation has $2n + \prod_{i \in [m]} (|b_i| + 1)$ variables and $3 + n + \newv{2}\prod_{i \in [m]} \newv{\frac{1}{2}}(|b_i| + 1)(|b_i| + 2) - 2\prod_{i \in [m]} (|b_i| + 1)$ constraints. \newv{Schaefer's \cite{Schaefer2009} formulation requires the constraints to be in inequality form $Ax \leq b, x \in \mathbb{Z}^n_+$.}

\newv{Our $\eqref{eq:igcr_lp}$ and Schaefer's \cite{Schaefer2009} formulations are larger than existing LP formulations for $\invlp$ \cite{Ahuja2001, Tavaslioglu2018}. For example, Ahuja and Orlin's \cite{Ahuja2001} LP formulation for $\invlp$ has at most $m + n$ constraints and $n$ variables. There exist polynomially solvable algorithms for $\invlp$ \cite{Ahuja2001}.}

\newv{From Table \ref{tbl:compare_lp}, our LP formulation for $\invgcr$ in \eqref{eq:igcr_lp} has many magnitudes fewer variables and constraints compared to Schaefer's \cite{Schaefer2009} formulation for $\invgcr$. By utilizing specific properties of $\gcrb$, our approach yields significantly smaller LP formulations for $\invgcr$ than those obtained by treating $\invgcr$ as an instance of a general inverse IP.}

\section{Conclusion}
\newv{Given a feasible basis of the LP relaxation of an IP, the GCR of the IP relaxes the nonnegativity constraints of the variables in the basis.  Then, the inverse GCR finds an objective vector (with minimum distance from a target objective vector) under which a given feasible solution $x^\circ$ for the forward GCR is optimal. If $x^\circ$ is also feasible for the IP, then we showed that solving the inverse GCR problems for all the feasible bases of the LP relaxation finds an upper bound on the optimal value of the inverse IP that is at least as tight as the optimal value of the inverse LP relaxation. We introduced cases where solving the inverse GCR problems for all the bases will exactly solve the inverse IP. If there exists a feasible basis of the LP relaxation where the components of $x^\circ$ in that basis are positive, then solving the inverse GCR defined by that basis will find an upper bound on the optimal value of the inverse IP that is at least as tight as the optimal value of the inverse LP relaxation. By reformulating the inverse GCR (defined on any basis) as an inverse shortest path problem, we presented an LP formulation for solving the inverse GCR. Our LP formulation for the inverse GCR has significantly fewer variables and constraints compared to an LP formulation for the inverse GCR obtained by applying inverse IP techniques proposed by Schaefer \cite{Schaefer2009}.






}



\section*{Acknowledgement}
\noindent This material is supported by the Office of Naval Research Grant N000142112262.

\printbibliography

\end{document}